\documentclass[amsa,nonatbib]{ipart}

\usepackage{graphicx}
\usepackage{hyperref}
\usepackage[table]{xcolor}
\usepackage[numbers,square]{natbib}
\usepackage{listings}
\startlocaldefs
\newtheorem{example}{Example}%

\endlocaldefs

\pubyear{2022}
\volume{0}
\issue{0}
\firstpage{1}
\lastpage{1}
\begin{document}

\begin{frontmatter}

%%  "Title of the Paper"
\title[Solving Nurse Scheduling Problem via PyQUBO]{Solving Nurse Scheduling Problem via PyQUBO
%\protect\thanksref{T1}
}
%\thankstext{T1}{Solving NSP via PyQUBO}

\begin{aug}
%%  \author{\fnms{John} \snm{Smith}\thanksref{t2}\ead[label=e1]{smith@foo.com}\ead[label=e2,url]{www.foo.com}}
%%  \thankstext{t2}{The author is supported by ...}
%%  \address{line 1\\ line 2\\ \printead{e1}\\\printead{e2}}

\author{\fnms{Matthew M.} \snm{Lin} \thanksref{t1} \ead[label=e1]{mhlin@mail.ncku.edu.tw}}
\address{Department of Mathematics, National Cheng Kung University,\\ Tainan 701, Taiwan\\\printead{e1}}
\author{\fnms{Yu-Chen} \snm{Shu}\thanksref{t2}\ead[label=e2]{ycshu@mail.ncku.edu.tw}}
\address{Department of Mathematics, National Cheng Kung University,\\ Tainan 701, Taiwan\\\printead{e2}}
\author{\fnms{Bing-Ze} \snm{Lu}\thanksref{t3}
        \ead[label=e3]{bingzelu@utexas.edu}}
\address{Oden Institute for Computational Engineering and
Sciences, The University of Texas at Austin,\\ TX 78712, USA\\\printead{e3}}
\and
\author{\fnms{Pei-Shan} \snm{Fang}\thanksref{t4}\ead[label=e4]{z11107031@ncku.edu.tw}}
\address{Department of Mathematics, National Cheng Kung University,\\ Tainan 701, Taiwan\\\printead{e4}}
\thankstext{t1}{The author is supported in part by the National Center for Theoretical Sciences of Taiwan and   
by the National Science and Technology Council of Taiwan
under grants 112-2636-M-006-002 and 112-2119-M-006-004.}
\thankstext{t2}{The author is supported in part by the National Science and Technology Council of Taiwan
under grants 112-2119-M-006-004 and 112-2115-M-006-001-.}
\thankstext{t3}{The author is supported in part by the National Science and Technology Council of Taiwan
under grant 113-2917-I-564-033.}
\thankstext{t4}{The author is supported in part by the National Science and Technology Council of Taiwan
under grants 112-2636-M-006-002 and 112-2119-M-006-004.}
\end{aug}
%%  History:
%\received{\sday{3} \smonth{1} \syear{2022}}

\begin{abstract}

The nurse scheduling problem is a critical optimization challenge in healthcare management. It aims to balance staffing demands, nurse satisfaction, and patient care quality. Corresponding to the constraints inherent in this scheduling problem, we detail the mathematical formulation step-by-step. We then utilize a quantum-inspired technique, the simulated annealing algorithm, and a quadratic unconstrained binary optimization model to optimize workload and increase nurse preferences. Numerical experiments are implemented to show the capacity of our proposed techniques. Our findings indicate a promising direction for future research, with  potential applications extending beyond nurse scheduling to other complex optimization problems.

%The nurse scheduling problem is a combinatorial optimization problem and can be solved using mathematical programming or heuristic approaches. It aims to find an optimal schedule that satisfies all the restrictions and requirements while optimizing a given objective function. 
%In this paper, we present a flexible quadratic unconstrained binary model with penalties to solve the problem effectively via PyQUBO.
%The objective function is a balanced combination of costs, nurse satisfaction, and workload. On the other hand, the solving progress can be accelerated by quantum-inspired techniques, such as digital annealing units. 
%We in this paper present a basic but flexible mathematical model to express the nurse scheduling problem. We show how non-constrained quadratic binary optimization can solve this problem effectively.
% We domonstrate the effectivity of using non-constrained quadratic binary method to solve the optimization problem
%In particular, we implement PyQUBO, an open-source Python library, to create the corresponding model and solve it via the simulated annealing algorithm. 
%Numerical experiments on real-world applications are implemented to show the capacity of our proposed techniques.
\end{abstract}


%%  Upper case for every keyword

%\tableofcontents

\end{frontmatter}

%%  The body
\section{Introduction}\label{sec1}

The Nurse Scheduling Problem (NSP) is NP-hard~\cite{OI00classification}. Finding an appropriate schedule can improve patient care and nursing staff satisfaction. A well-designed schedule should ensure that enough nurses are present to meet the patients' needs, without sacrificing the nurses' availability and preferences. In contrast, inadequate scheduling can result in a heavy workload, extended shifts, burnout among nurses, and negative outcomes in patient care. Thus, developing efficient tools to solve nurse scheduling problems is important not only in increasing the nurses' well-being but also in the quality of care for every patient. 
There are complex constraints involved in scheduling shifts for nurses. For example, in~\cite{STB13generic} these constraints include the number of shifts per nurse, minimum staffing levels for multiple shifts or different wards, the demand for specific nursing skills, and nurse preferences. 
Although the NSP is an NP-hard problem, one can apply 
 algorithms and heuristic methods  to develop near-optimal solutions. 
For example, in~\cite{Legrain201}, a simple heuristic approach is presented and can be easily implemented by operating on spreadsheets without the need for expensive commercial software. 
In~\cite{WHZO14solving}, a binary goal programming model using optimization software LINGO~\cite{Lingo2002} is presented to study the NSP in outpatient departments. This model categorizes nurses as either full-time or part-time, with different limitations, while considering nurses' exceptional skills and diverse tasks. In~\cite{AGF17new}, a multi-commodity network flow model with distinct sources (nurses) and sink nodes (shifts) is given in a real case study of the NSP in an Egyptian hospital. This model is verified to obtain a schedule with the overtime cost per week 36\% less than the one manually created by the supervisor head nurse. In~\cite{GHG20simple}, a generic variable-fixing method is  applied to reduce the size of the NSP and make it more comfortable to solve by discarding nurses, shift patterns, and binary variables. After that, a general-purpose mixed-integer programming solver is utilized to solve this reduced problem. As noted in~\cite{GHG20simple}, this reduced model offers an alternative approach to approximating the optimal solution for the original NSP, but it does not provide a mathematically equivalent optimal solution to the original. %See also~\cite{BCBL04state, STB13generic} for more information about the NSP.

In this work, we consider the NSP as a constrained quadratic unconstrained binary optimization (QUBO) problem. We showcase the employment of "PyQUBO," an open-source Python library for the QUBO, with an implementation of the simulated annealing algorithm (SAA) to find an optimal solution that satisfies hospital policies and nurses' preferences.
Our work includes two steps.
The first step is a problem-defining process. We use binary variables, where the value 1 represents a nurse assigned to a shift and the value 0 indicates the other way around, to describe the nurse's assignment. We represent the requirements, including the number of nurses for each shift, the number of consecutive shifts a nurse can accept, and the number of days off required by policies, as equality and inequality constraints, notably a series of logic constraints. 
The second step is to formulate the problem in PyQUBO by using the PyQUBO library to define the variables and the objective function in PyQUBO's syntax. Since the NSP problem is a QUBO problem, we have to convert the constraints into the objective function by adding the product of Lagrange multipliers and constraints to the objective function. Once set up the QUBO problem, we minimize it with the implementation of the SAA, which is a classical optimization solver for solving QUBO problems.

In Taiwan, there are regulations regarding work shifts, which specify the order of consecutive workdays and the mandatory rest periods.To our knowledge, no existing research applies the QUBO framework and the simulated annealing algorithm, let alone quantum annealing, to create hospital shift schedules. We aim to increase nurse scheduling flexibility in Taiwan, empowering nurses to incorporate their preferences into scheduling rules more easily. The QUBO structure enables the conversion of rules into mathematical expressions, and the simulated annealing method efficiently identifies local minima, making it an ideal choice. We utilize these tools to alleviate the workload of individuals responsible for meeting all obligations while also ensuring the smooth functioning of hospitals and complying with labor laws. There is yet no standard yardstick available for measuring the effectiveness of traditional computing and algorithms as compared to quantum computing in addressing QUBO problems. We believe that this research will provide valuable perspectives.

The paper is organized as follows. In Section 2, we provide an overview of the fundamental concepts related to the optimization problem, the annealing algorithm, and the Python package used for handling the QUBO problem. In Section 3, we describe our specific problem using concrete mathematical language. In Section 4, we outline the step-by-step process of employing PyQUBO to solve the NSP with a detailed example. Section 5 demonstrates the application of the optimization method, adhering to the specified rules, to formulate a suitable nurse schedule. Lastly, in Section 6, we conclude our findings and discuss potential future extensions of this research. 
 
\section{Preliminaries}\label{sec2}
We introduce the definition of {the} QUBO, SAA, and the state-of-the-art Python library, PyQUBO. The understanding of these tools can facilitate our later discussion.
\subsection{Quadratic Unconstrained Binary Optimization Problem}\label{subsec2.1}
The QUBO problem is a combinatorial optimization problem for solving the following optimization problem involving a quadratic function and binary variables.
 
\begin{equation}\label{eq:QUBO}
\begin{array}{cl}
	\mbox{Minimize} & f(\mathbf{x}):=\mathbf{x}^\top Q \mathbf{x} 
=   \sum\limits_{i, j} q_{ij} x_ix_j,\\
       \mbox{subject to}& 
       \mathbf{x}\in \{0,1\}^n,
	\end{array}
\end{equation}
where $Q = [q_{ij}]\in \mathbb {R} ^{n\times n}$ is a %n upper triangular
matrix with entries $q_{ij}$ indicating a weight for each pair of indices $i,j\in \lbrace 1,\dots ,n\rbrace$ and $\mathbf{x} = [x_i] \in \{0,1\}^n$.
This concern has a wide range of applications, including the traveling salesman problem~\cite{papalitsas2019qubo}, the knapsack problem~\cite{GKAA02solving}, and the graph coloring problem~\cite{kochenberger2005unconstrained}. It is beneficial for operations research tasks, such as scheduling and resource allocation~\cite{GE18scheduling}. For a comprehensive overview of the theoretical results and various applications of QUBO problems, one can refer to~\cite{MANY14unconstrained}.

In the QUBO framework, the problem can be interpreted as a weighted graph $G = (V, E)$, where $V = \{v_1, \ldots, v_n\}$ denotes the set of vertices and $E = \{(v_i, v_j)\}_{i, j = 1}^n$ symbolizes the edges connecting vertices $v_i$ and $v_j$, with $q_{ij}$ denoting the weight on edge $(v_i, v_j)$. If $q_{ij} = 0$, it indicates the absence of the edge $(v_i, v_j)$ in $E$. The QUBO framework aims to find a subset $S \subset V$ that minimizes the sum of the weights in the induced subgraph. The indicator variable $x_i = 1$ signifies the inclusion of vertex $v_i$ in subset $S$, and $x_i = 0$ otherwise. The QUBO formulation has gained attention due to its relation to the Ising model, a fundamental model in statistical mechanics that describes the energy interactions in a system of atoms with positive or negative spins. The Ising model is expressed mathematically as:
\begin{equation*}
    E = -\sum_{i, j} J_{i,j} s_i s_j - \sum_i h_i s_i
\end{equation*} 
where $E$ represents the system's energy, $s_i$ denotes the spin of the $i$-th atom taking values in $\{-1, 1\}$, $J_{ij}$ is the interaction energy between spins, and $h_i$ is the external magnetic field influencing each spin. The summation over indices $i, j$ runs over neighboring pairs of atoms. The Ising model can be reformulated as:
\begin{equation*}
E = \mathbf{s}^T J \mathbf{s}, 
\end{equation*}
where $\mathbf{s}$ is a vector composed of spin values $s_1, \ldots, s_n$, and the interaction weights $J_{i, j}$ are elements of the matrix $J$ at positions $(i, j)$. Furthermore, the external magnetic field effects, represented by $h_i$, are integrated into the diagonal elements of the matrix $J$. This reformulation enables a compact matrix representation of the Ising model and simplifies analysis and computation.

The Ising model, like the QUBO, can also be represented as a graph. The goal is to determine a spin configuration that minimizes the overall energy of the system. This graph-based approach is similar to that of the QUBO and emphasizes the structural similarities between the two models. The only difference is that the QUBO uses binary variables $x_i \in \{0, 1\}$ while the Ising model uses spin variables $s_i \in \{-1, 1\}$. This conversion can be achieved through the relation 
\begin{equation*}
    s_i = 2x_i - 1.
\end{equation*} 
Substituting $x_i = \frac{1}{2}(s_i + 1)$ into the QUBO objective function $x^\top Q x$ and rearranging terms allows the QUBO problem to be represented in the form of the Ising model $E = \mathbf{s}^T J \mathbf{s}$, where the elements of matrix $J$ and the vector $\mathbf{s}$ are adjusted accordingly to reflect this transformation. 
The relationship between the QUBO and Ising model demonstrates the adaptability of these frameworks for representing combinatorial optimization problems. On the other hand, the Ising model can be effectively solved using quantum algorithms, which offer a quick solution for certain cases. This connection opens up a promising way of solving the QUBO problem by converting them into the Ising model and making use of quantum computational techniques. 
  
So far, quantum computing has the potential to greatly reduce computation time. However, it is important to recognize the limitations of quantum computers, which give challenges in controlling noise and errors, as well as the loss of quantum properties during calculations and high energy consumption. This study focuses mainly on utilizing the SAA to solve the QUBO problem, instead of directly developing a quantum algorithm. This approach was chosen to effectively address the problem while mitigating the challenges associated with quantum computing technologies.

\subsection{Simulated Annealing Algorithm}\label{subsec2.2}
Developed in 1983 by Kirkpatrick, Gelatt, and Vecchi~\cite{KGB83optimization},  the SAA is applied to solve the traveling salesman problem by mimicking the process of heating and cooling a physical system, where the temperature progressively decreases over time.
% 我覺得可以省略 from an initial positive value. 
The algorithm initiates with a state $s_0$ with a high temperature to enlarge the scope of the search space and then gradually reduces the temperature by randomly jumping to %its
the neighboring state $s_1$ of $s_0$. If the nearby state $s_1$ has a lower objective value, the algorithm replaces $s_0$ with $s_1$
and repeats the iteration from $s_1$; otherwise, it reselects another nearby state as $s_1$ and repeats the discussion stated above. However, the SAA also picks up a nearby state $s_1$ to replace $s_0$ concerning the determination of an acceptance probability function even if the computed value is larger at state $s_1$. This feature can prevent the SAA from terminating prematurely at a local minimum of the objective function. This above process will terminate until the number of iterations is reached.% 

In addition, the convergence of the SAA is guaranteed in \cite{VA87simulated}. %
The study shows the convergence of SAA is four times faster than that of tabu search~\cite{Beasley98heuristic}.
In~\cite{KN01performance}, the SAA is applied to some benchmark QUBO problems. These experiments show that the SAA has better performance and shorter running time than those searching methods.
These advantages then motivate us to use
the SAA for solving the NSP in this work.

\subsection{PyQUBO}\label{subsec2.3}

This work employs PyQUBO, a Python library, to solve QUBO problems. PyQUBO is recognized for its intuitive and user-friendly programming interface, which facilitates the straightforward definition of binary variables within mathematical models. It incorporates the following key features.

First, PyQUBO supports two distinct classes of representations: the spin direction $\{-1, 1\}$ for the Ising model and the binary variable $\{0, 1\}$ for QUBO problems. This versatility facilitates the easy conversion between these data structures, allowing users to define variables as spin states or binary variables effortlessly.
%
% Second, 
After setting the variables, it provides a `.compile()` function that translates the defined objective function into a structured data format, which can then be transformed into a QUBO or an Ising model. This process expands the objective function into a weighted graph format, where the nodes represent the variables, and the edges are weighted by the coefficients of variable pairs to define their interactions. Once this step is completed, the model formulation process is finished.

% {\color{red}Second, it integrates with D-Wave's "neal" package to facilitate the sampler of the SAA for solving QUBO problems. 
% This integration incorporates a sampling function for initializing solution points and a updating procedure to find an approximate solution to the optimization problem. A simple code is listed below. The object function is given in the variable `ObjectFuncton` in prior and be omitted here.
Second, PyQUBO integrates with D-Wave's "neal" package to facilitate the SAA sampler for solving QUBO problems. 
This integration includes a sampling function that helps initialize solution points and an updating procedure to find an approximate solution to the optimization problem. A simple code is listed below. The object function is given in the variable `ObjectFuncton` in the prior and is omitted here.
\begin{lstlisting}[frame=single]
import neal   
sampler = neal.SimulatedAnnealingSampler() # create the SAA sampler using the 'neal' package
model = ObjectFunction.compile()
sampleset = sampler.sample(model.to_bqm(), num_reads=1000)
decoded_samples = model.decode_sampleset(sampleset)
best_sample = min(decoded_samples, key=lambda x: x.energy)
\end{lstlisting}
After completing the necessary preparations, the `min` function is used to find the minimum value of the objective function over a specified number of samplings.
% Following these preparations, the `min` function at the end is to find the minimum of the designed objective function over a specific number of sampling. 
% }

Third, PyQUBO is a powerful tool that allows users to translate logic operations into binary representations. This feature is handy for optimization problems that involve multiple constraints and logical conditions. For example, exclusive conditions are necessary to correctly model the problem in scheduling scenarios where nurses should not operate consecutively. In other words, it provides functions like NotConst, OrConst, AndConst, and XorConst to help formulate problems requiring logical representation. 

The following section will mention how to modify the constraint problem to an unconstrained one. These functions enable users to directly integrate complex logical constraints into the binary optimization framework, thus enhancing the expressiveness and applicability of the QUBO formulations. We will utilize the previously mentioned functions to address the nurse scheduling challenge, wherein binary variables represent the presence of a nurse on duty. These values of the variables should align with the hospital's operational guidelines and facilitate the creation of a comprehensive shift schedule, serving as a valuable reference for coordinating the nurses' rotations.

\section{Problem Description and Mathematical Model}\label{sec3}
In Taiwan, a three-shift system, day shift, night shift, and graveyard shift (late night shift), is widely used for nurse scheduling. Suppose that there are $N$ nurses for arrangement.
%in a month 
% and $n$ nurses on duty per day. 
 Our goal is to maximize the number of consecutive leaves for each nurse and to create 
 % a proper
an appropriate shift table that does not fail the following rules:
%\textbf{Rules:}
\begin{enumerate}
    %\item Each nurse has at most one shift per day.
    \item [1.] $m_1$ nurses are assigned to work graveyard shifts %{\color{red}per month} 
    with $n_1$ of them on duty per day, where $0 \leq n_1 < m_1$.
    \item [2.]$m_2$ nurses are assigned to work night shifts with $n_2$ of them on duty per day, where $0 \leq n_2 < m_2$.
    \item[3.] $n_3$ nurses will be on duty for day shifts per day.
    \item[4.] Nurses who work graveyard or night shifts should work at least two-day shifts, but no more than $k$-day shifts, where $k \geq 2$, for each arrangement. After 2 to $k$ consecutive shifts, each nurse must have a leave.
    \item [5.] Nurses who work on the weekends can take a break for at least two days during the week.
    % , and there are no night and graveyard shifts between Saturday and Sunday evenings. 
    Here, we assume that the first day of the week is Saturday and the last day of the week is Friday.
\end{enumerate}

Let $d$ represent the total number of days required for the arrangement.
Facing rules 1-5, we can express the nurse scheduling problem as a constrained combinatorial optimization problem in the form of 
\begin{equation}\label{eq:combinatorial}
\begin{array}{cl}
	\mbox{Minimize} & H_{{cost}}(\mathbf{x})\\
       \mbox{subject to}& H_{{const}}(\mathbf{x}) \mbox{ and }
       \mathbf{x}\in \{0,1\}^{N\times d},
	\end{array}
\end{equation}
where the objective function $H_{{cost}}(\mathbf{x})$ represents one type of "soft" constraint in the scheduling problem, and the constraint $H_{{const}}(\mathbf{x})$ is defined according to rules 1-5. The fundamental problem of interest is to assume $H_{cost}(\mathbf{x}) = 0$ and focus on satisfying the constraint $H_{{const}}(\mathbf{x})$. 
Note that~\eqref{eq:combinatorial} is a constrained optimization problem, while a QUBO model in~\eqref{eq:QUBO} is to solve a binary combinatorial optimization problem without any constraints. 
To utilize QUBO to solve~\eqref{eq:combinatorial}, we consider minimizing the objective function  
\begin{equation}\displaystyle
\label{one const}
    H(\mathbf{x}) = H_{{cost}}(\mathbf{x}) +  \lambda^\top H_{const}(\mathbf{x})
\end{equation}
over $\mathbf{x} \in \{0,1\}^{N\times d}$ while given with a suitably selected constraint term weight $\lambda$.

Certainly, $ 0 < n_1 + n_2 + n_3 < N$. To reduce the number of variables, we start our discussion by dividing the nurses into three groups, including 
the graveyard shift group $U_1$, the night shift group $U_2$, and the day shift group $U_3$. We use the notation $p_i$ for $i=1,\ldots, N$ to represent each nurse.
Without loss of generality, we let the group $U_1$ consist of nurses $p_i$, $i=1,\ldots,m_1$, the group $U_2$ consist of nurses $p_{i}$, $i=m_1+1,\ldots,m_1+m_2$, and the group $U_3$ 
consist of the rest of the nurses. 
% consisting of nurses $p_{i}$, $i=m_1+m_2+1,\ldots,N$. 
Let $x_{ij} \in \{0,1\}$ be the variable that represents whether the nurse $p_i$ will work on day $j$ or not, that is,
\begin{equation*}
    x_{ij} = \left\{\begin{array}{ll}
                 1, & \mbox{if nurse $p_i$ is on duty on day $j$} \\  
                 0,      & \mbox{otherwise.}
                \end{array} \right.
\end{equation*}
Since the main concern is the determination of $H_{const}$, we carry out such a task by $H_{const}$ in four phases: $T_i$, $i=1,\ldots,4$,
(i.e., $H_{const} = [T_1,T_2,T_3,T_4]^\top$).

Note that the number of nurses in each group is subject to the limitations given by 
% the rules: 1 to 3 
rule 1 to rule 3 and can be illustrated by the following 
% three types of 
constraints:
%\begin{subequations}
%\label{rule234}
\begin{align}\displaystyle
\left\{
    \begin{array}{ll}
       h_j^{(1)}(\mathbf{x}) :=  \sum\limits_{i = 1}^{m_1} x_{ij} - n_1 = 0, & j = 1,\ldots,d,\\
       h_j^{(2)}(\mathbf{x}) := \sum\limits_{i = m_1+1}^{m_1+m_2} x_{ij} - n_2 = 0,& j = 1,\ldots,d,\\
       h_j^{(3)} (\mathbf{x}):= \sum\limits_{i = m_1+m_2+1}^{N} x_{ij} - n_3 = 0,&  j = 1,\ldots,d,
    \end{array}
    \right. \label{rule234}
    \end{align}
    with $\mathbf{x} = [x_{i,j}]$.
%\end{subequations}
% Putting 
Arraying all these elements together, we then define the following constrained formulation:
\begin{equation} 
\label{T1const}
T_1(\mathbf{x}) = \sum\limits_{j = 1}^{d} \left ( \left[h_j^{(1)}(\mathbf{x})\right]^2+ \left [h_j^{(2)}(\mathbf{x})\right ]^2 + \left[h_j^{(3)}(\mathbf{x})\right]^2\right )
\end{equation}

In the second phase, let $\mbox{AND}$, $\mbox{OR}$, and $\mbox{NOT}$ be the logical gates. 
% The 
These analytical representations of logic gates are
% can be 
defined by
 \begin{subequations}
     \label{eq:logic}
\begin{align}\displaystyle
%\left\{
&\mbox{AND}(a,b) = ab, \\
&\mbox{OR}(a,b) = a+b-ab, \\
&\mbox{NOT}(a) = 1-a, 
 %   \right. 
    \end{align}
    \end{subequations}
where $a$ and $b$ are two binary inputs.\\   
In PyQUBO, the logic gates are expressed by three class inputs (or two in the case of NOT), i.e., $\mbox{AndConst}(a,b,c,\mbox{\textquotesingle and\textquotesingle})$, $\mbox{OrConst}(a,b,c,\mbox{\textquotesingle or\textquotesingle})$, or $\mbox{NotConst}(a,b,\mbox{\textquotesingle not\textquotesingle})$, respectively,
where $c$ is 
% another 
an auxiliary binary input. When the binary variables satisfy the constraints, 
$\mbox{AND}(a,b)= c$, $\mbox{OR}(a,b)= c$, or $\mbox{NOT}(a)= b$, the resulting energy in the annealing algorithm is 0; 
otherwise
% Otherwise, 
the output energy is 1 \cite{ZTT21pyqubo}.  
%這一句我有點不太懂意思
Concerning rule 4, we divide this constraint into two parts. The first part deals with the constraint of "working at least two-day shifts." 
The constraint requires composited strings of $x_{ij}$'s for each $i$  in such a way that none of them includes $10$ in the first two entries and $01$ in the last two entries, and they are not allowed to have $010$ as a substring.
% neither includes $10$ in the first two entries and $01$ in the last two entries nor permits $010$ to be its substring.
% This means that for each $i$, the string for variables $x_{ij}$'s does not include the condition that the first and last two entries have entries $01$ and $10$, respectively, or a substring has entries equal to $010$.
Using the logic gate classes, we can define the logical constraints for the cases 010, 10, and 01 %, respectively,
%working at least two-day shifts 
in the following logical constraints: 
 \begin{subequations}
\label{rule5}
    \begin{align}\displaystyle
   &     \mbox{AndConst}\left(x_{i,j+1}, 1 - \mbox{OrConst}(x_{ij},x_{i,j+2},1,  \mbox{\textquotesingle or\textquotesingle}), 0, \mbox{\textquotesingle and\textquotesingle}\right), \\
  &      \mbox{AndConst}\left(x_{i,1}, \mbox{NotConst}(x_{i,2},1,\mbox{\textquotesingle not\textquotesingle}), 0,  \mbox{\textquotesingle and\textquotesingle}\right), \\
  &      \mbox{AndConst}\left(x_{i,d}, \mbox{NotConst}(x_{i,d-1},1,\mbox{\textquotesingle not\textquotesingle}), 0,\mbox{\textquotesingle and\textquotesingle}\right), 
    \end{align}
\end{subequations}
for $i \in \left\{1, ..., m_1+m_2\right\}$ and $j \in \left\{1,...,d-2 \right\}$. 
Without causing any ambiguity, we rewrite~\eqref{rule5} by using the logic gates given in~\eqref{eq:logic} as
    \begin{align*}%\displaystyle
&     \ell_{i,j}(\mathbf{x}):=   \mbox{AND}\left(x_{i,j+1}, 1 - \mbox{OR}(x_{i,j},x_{i,j+2})\right) = 0,\\
&      \ell_{i,1}(\mathbf{x}):=  \mbox{AND}\left(x_{i,1}, \mbox{NOT}(x_{i,2})\right) = 0,\\
 &    \ell_{i,d}(\mathbf{x}):= \mbox{AND}\left(x_{i,d}, \mbox{NOT}(x_{i,d-1})\right) = 0,
    \end{align*}
for $i \in \left\{1, ..., m_1+m_2\right\}$ and $j \in \left\{1,...,d-2 \right\}$, respectively. Then, the constraint $T_2$ is given by
\begin{equation*}
T_2(\mathbf{x}) := \sum_{i = 1}^{m_1+m_2} \sum_{j = 1}^{d}\ell_{i,j}(\mathbf{x}).
\end{equation*}

We have already discussed how to incorporate equality constraints into the objective function.
In the third phase, we handle the constraint "no more than $k$-day shifts." This statement indicates that nurses must have at least one day off after working for $k$ consecutive days and can be expressed in terms of the following constraint:
\begin{align}\displaystyle
\label{nomorekdays}
    \sum_{j = \ell}^{\ell+k+1} x_{ij} \leq k, \quad \mbox{for } i= 1, \ldots, m_1+m_2\mbox{ and } \ell = 1,\ldots,d-(k+1).
\end{align}
 Note that~\eqref{nomorekdays} represents a set of inequality constraints. The direct way of converting an inequality constraint to an equality one is to add additional slack variables. For instance, consider an inequality constraint:
\begin{equation}
\label{ineq}
    g(\mathbf{x}) \leq {b},
\end{equation} 
for some nonnegative integer $b$. The inequality~\eqref{ineq} implies that the values of $g(\mathbf{x})$ must 
% be lying
lie in the range of $[0,b]$. By adding the slack variables, we 
% can 
then transform the representation of~\eqref{ineq} into the equality constraint:
\begin{equation*}
\label{slack2}
    g(\mathbf{x}) =  s_1 + ... + s_b,
\end{equation*}
where $s_i \in \{0,1\}$ for $i=1,\ldots,b$. 

% Assuming $g_{i\ell}(\mathbf{x})=\sum_{j = \ell}^{\ell+k+1} x_{ij}$, $i= 1, \ldots, m_1+m_2$ and $\ell = 1,\ldots,d-(k+1)$, we l
Let 
$s_p^{i,\ell}$, $p = 1,\ldots,k$, be the corresponding $k$ independent slack variables for each summation in~\eqref{nomorekdays}.
Then for  $i= 1, \ldots, m_1+m_2$ and  $\ell = 1,\ldots,d-(k+1)$, we can represent~\eqref{nomorekdays} in terms of equality constraints denoted by
\begin{align*}\displaystyle
\label{nomorekdays2}
    g_{i\ell}(\mathbf{x};\mathbf{s}^{i,\ell}) := 
    \sum_{j = \ell}^{\ell+k+1} x_{ij} - 
    \sum_{p=1}^ks_p^{i,\ell} = 0, 
\end{align*}
 where the vector $\mathbf{s}^{i,\ell} = [s_p^{i,\ell}]$ and 
 each entry $s_p^{i,\ell} \in \{0,1\}$. With the above substitution, we can encode the condition $T_3$ in terms of
\begin{equation*}
T_3 (\mathbf{x};\mathbf{s}) = \sum_{i=1}^{m_1+m_2} \sum_{l = 1}^{d-(k+1)} \left(g_{i\ell}(\mathbf{x};\mathbf{s}^{i,\ell})\right)^2,
\end{equation*}
 where $\mathbf{s} = \left[\mathbf{s}^{1,1},\ldots,\mathbf{s}^{1,d-(k+1)},\mathbf{s}^{2,1},\ldots, \mathbf{s}^{m_1+m_2,d-(k+1)} \right]$. 
% In the forth phase, we note that for rule 5, we restrict nurses in group $U_3$ by a stronger constraint, ``Nurses working on the weekends can take a break for at least two days on weekdays". 
% This means that if nurses have to work on Saturday or Sunday, we have to assign at least two-day breaks to them this week.
% Let $U_{Sat}$ be an index set tagging Saturdays in a month. 

In the fourth phase, we see that for rule 5, if nurses work on Saturday or Sunday, we have to assign at least two-day breaks to them in the same week. Let $U_{Sat}$ collect the index of Saturdays in a month. We then have the following constraints:  
\begin{align}
\label{rule6}
    \sum_{j = l}^{l+6} x_{ij} \leq 5,  \quad \mbox{for } i \in \{1, ... , n\} \mbox{ and } l \in U_{Sat},
\end{align}
to define rule 5. For this matter, suppose $\ell$ and $i$ are fixed; we can identify~\eqref{rule6} as the equality constraint given by
\begin{align*}\displaystyle
\label{nomorekdays3}
    h_{i\ell}(\mathbf{x};\mathbf{y}^{i,\ell}) := 
    \sum_{j = \ell}^{\ell+6} x_{ij} - 
    \sum_{p=1}^5 y_p^{i,\ell} = 0, 
\end{align*} 
by adding five slack variables $y_p^{i,\ell}\in \{0,1\}$, $p=1,\ldots,5$, with $\mathbf{y}^{i,\ell} = [y_p^{i,\ell}]$. The restriction of $T_4$ is given by 
\begin{equation*}
T_4(\mathbf{x};\mathbf{y}) = \sum_{i=1}^{n} \sum_{\ell  \in  U_{Sat}}
\left(h_{i\ell}(\mathbf{x};\mathbf{y}^{i,\ell})\right)^2.
\end{equation*}
 where $\mathbf{y}$ represents a collection of variables $\mathbf{y}^{i,\ell}$, $i=1,\ldots, n$ and $\ell  \in  U_{Sat}$.

\section{A PyQUBO Model for NSPs}\label{sec4}
The following is an example of a PyQUBO formulation for the nurse scheduling problem. We start our illustration with the assignment of five nurses working the graveyard shift for five days. We assume that each day only three nurses are required to work. It follows from a direct computation that the best schedule is to assign each nurse to work for three days in every five-day arrangement. Below, we show how the NSP problem can be written and solved using PyQUBO.

We start with importing the needed packages and inputting the given parameters.\\
%\newpage
\lstset{frame=single,framexleftmargin=-1pt,framexrightmargin=-17pt,framesep=12pt,linewidth=0.98\textwidth,language=pascal}
\begin{lstlisting}[frame=single]
## import packages
from pyqubo import Binary, Constraint, Array
import numpy as np
import neal  # the SAA sampler
sampler = neal.SimulatedAnnealingSampler()
## input parameters
m1 = 5 # the number of nurses for arrangement
n1 = 3 # the number of nurses required per day 
d = 5 # the number of days required to be determined
\end{lstlisting}

After setting up the packages and parameters, we define a binary variable, X, to record whether the nurses are on duty. For convenience, we reshape the array into a $m_1\times d$ matrix, where the $i$th row represents the $i$th nurse, and the $j$th column represents the $j$th day. We also use the binary variables, total\_shift and sum\_shift, to collect the number of nurses per shift and the number of shifts arranged for nurses.

\begin{lstlisting}[frame=single]
## X: binary variables for saving nurses' allocation
X_initial = Array.create("Graveyard", shape = m1*d, vartype = "BINARY")
X = np.zeros(m1*d).reshape(m1, d)
X = X.tolist()
for i in range(m1):
    for j in range(d):
        X[i][j] = X_initial[d*i+j]    
## total_shift: number of nurses per shift
total_shift = []
for i in range(d):
    col_sum = 0
    for j in range(m1):
        col_sum = col_sum + X[j][i]
    total_shift.append(col_sum)
## sum_shift: number of shifts arranged for nurses
sum_shift = []
for i in range(m1):
  sum_shift.append(sum(X[i][j] for j in range(d)))
\end{lstlisting}

We use the following constraint to represent our $T_1$ constraint defined in~\eqref{T1const}.
\begin{lstlisting}[frame=single]
## the limitation of manpower
Ha = sum((total_shift[i] - n1)**2 for i in range(d))
\end{lstlisting}
We also add one more feature to our demonstration to demonstrate the potential of our code. To this, we include the need for each nurse to work for only three days in this five-day arrangement into our constraint as below.
\begin{lstlisting}[frame=single]
## the constraint for satisfying the optimal condition
Hb = sum((sum_shift[i] - 3)**2 for i in range(m1))
\end{lstlisting}
To transform the current constrained problem into an unconstrained one, we encapsulate the above two constraints to obtain and compile the objective function. The records for the nurse scheduling are generated by the SAA ten times, while 10000 sweeps are used to calculate the result each time.
% ycshu@0505: I modify the H to be ObjectFunction to be consistent with the previous introdction.
\begin{lstlisting}[frame=single]
## define and compile the objective function
ObjectFunction = Constraint(Ha,"eachshift") + Constraint(Hb,"eachworker") 
model = ObjectFunction.compile()
bqm = model.to_bqm() # bqm: binary quadratic model 
sampleset = sampler.sample(bqm, num_reads=10, num_sweeps=10000)
# num_reads: the number of experiments
# num_sweeps: the number of sweeps used in each calculation.
decoded_samples = model.decode_sampleset(sampleset)

## pick up the record with the lowest energy
graveyard_sampleset = min(decoded_samples, key=lambda s: s.energy)
graveyard_record = graveyard_sampleset.sample
## prepare the optimal scheduling table
graveyard_table = np.zeros(m1*d)
for key, value in graveyard_record.items():
    newkey = int(key.replace("Graveyard[", "").replace("]", ""))
    graveyard_table[newkey] = value
graveyard_table = graveyard_table.reshape(m1, d).astype(int)
print(graveyard_table)
\end{lstlisting}
Here, we finally use the~\emph{numpy} package to generate a binary matrix, graveyard\_table, with the value equal to 1, indicating that one nurse is on duty, and the value equal to 0, indicating that one nurse is off duty. Finally, we can check the satisfaction of the optimal condition by using the following command in PyQUBO. 
\begin{lstlisting}[frame=single] 
## Check whether the result fits the constraints.
print(graveyard_sampleset.constraints())       
\end{lstlisting}
If the conditions are met, the output values for the constraints will equal zero, which are just what we see in this simple example.

\section{Numerical Experiments}\label{sec5}
We have defined the NSP in Section 3, and we demonstrate how the NSP can be written and solved using the PyQUBO in Section 4. In this section, 
we want to analyze two patterns from Section 3 with two different limitations on the maximal number of consecutive working days, i.e., $k = 4$ or $5$. For different $k$, we repeat our experiments ten times with 10000 sweeps for each calculation. We choose the result with the lowest energy from these ten results. 
To ease our discussion, we collect related parameters from Section 3 in Tabel~\ref{tab:par}.

\begin{table}[ht!!]
       \caption{Parameters for describing mathematical models.}
    \centering
    \begin{tabular}{|l|m{8cm}|}
        \hline
constant parameters                        & description                                                                                                    \\ \hline\hline
        $N$                              & the total number of nurses
 % \\\hline
        % $n$                             & the required number of nurses on duty per day
                                   \\\hline
        $m_i$, $i=1,2$                              & the number of nurses assigned to work graveyard and night shifts, respectively       
                                        \\\hline
        $n_i$, $i=1,\ldots,3$                              & the number of nurses working graveyard, night, and day shifts should appear per day, respectively                                          \\\hline
        $k$                              & the number of consecutive working days                                                                                                                                      \\\hline
        $d$ & the number of days required for the arrangement
        \\\hline
    \end{tabular}
 \label{tab:par}
\end{table}

\begin{example}

% In this section, we show two patterns of the results, one with the parameter $k = 4$, and the other with $k = 5$. 

We consider the arrangement for September 2022, where $d = 30$ and the first day of this month is Thursday, and let the parameters $N = 13$, $m_1 = m_2 = 3$, $m_3 = 7$, $n_2 = n_2 = 2$, $n_3 = 3$, respectively. To present the results clearly, we split the scheduling tables into three subtables corresponding to  $U_1$, $U_2$, and $U_3$, respectively. In the tables, the gray block means that the nurse is on duty for that shift, and the white block means that the nurse is on leave. We use the PyQUBO, as demonstrated in Section 3, to express this QUBO with a simulated annealing-based sampler, Neal~\cite{neal2023}.

 \begin{itemize}
\item {Model 1: k = 4}
\begin{itemize}
    \item $U_1$: Graveyard-shift group\\[1em]
    \resizebox{.9\textwidth}{!}{%
    \begin{tabular}{|c|c|c|c|c|c|c|c|c|c|c|c|c|c|c|c|}
    \hline
       Date & 1& 2 & 3 & 4 & 5& 6& 7& 8 & 9 & 10 & 11 & 12 & 13 & 14 & 15\\
    \hline
       $P_1$ & \cellcolor{gray} & \cellcolor{gray} &   & \cellcolor{gray} & \cellcolor{gray} &   & \cellcolor{gray} & \cellcolor{gray} &  \cellcolor{gray} &   &   & \cellcolor{gray} & \cellcolor{gray} & \cellcolor{gray} & \cellcolor{gray}\\
    \hline
       $P_2$ & \cellcolor{gray} & \cellcolor{gray} & \cellcolor{gray} &   &   & \cellcolor{gray} & \cellcolor{gray} & \cellcolor{gray} &   & \cellcolor{gray} & \cellcolor{gray} &   & \cellcolor{gray} & \cellcolor{gray} &  \\
    \hline
       $P_3$ &   &   & \cellcolor{gray} & \cellcolor{gray} & \cellcolor{gray} & \cellcolor{gray} &   &   & \cellcolor{gray} & \cellcolor{gray} & \cellcolor{gray} & \cellcolor{gray} &   &   & \cellcolor{gray}\\
    \hline   
    \hline
       Date & 16 & 17 & 18 & 19 & 20 & 21 & 22 & 23 & 24 & 25 & 26 & 27 & 28 & 29 & 30\\
    \hline
       $P_1$ &   & \cellcolor{gray} & \cellcolor{gray} & \cellcolor{gray} & \cellcolor{gray} &   &   & \cellcolor{gray} & \cellcolor{gray} & \cellcolor{gray} &   & \cellcolor{gray} & \cellcolor{gray} &   &  \\
    \hline
       $P_2$ & \cellcolor{gray} & \cellcolor{gray} &   &   & \cellcolor{gray} & \cellcolor{gray} & \cellcolor{gray} &   &   & \cellcolor{gray} & \cellcolor{gray} &   & \cellcolor{gray} & \cellcolor{gray} & \cellcolor{gray}\\
    \hline
       $P_3$ & \cellcolor{gray} &   & \cellcolor{gray} & \cellcolor{gray} &    & \cellcolor{gray} & \cellcolor{gray} & \cellcolor{gray} & \cellcolor{gray} &    & \cellcolor{gray} & \cellcolor{gray} &   & \cellcolor{gray} & \cellcolor{gray}\\
    \hline   
    \end{tabular}}\\[1em]
    \item $U_2$: Night-shift group\\[1em]
    \resizebox{.9\textwidth}{!}{%
     \begin{tabular}{|c|c|c|c|c|c|c|c|c|c|c|c|c|c|c|c|}
    \hline
       Date & 1 & 2 & 3 & 4 & 5 & 6 & 7 & 8 & 9 & 10 & 11 & 12 & 13 & 14 & 15\\
    \hline
       $P_4$ &   &   & \cellcolor{gray} & \cellcolor{gray} & \cellcolor{gray} & \cellcolor{gray} &   &   & \cellcolor{gray} & \cellcolor{gray} & \cellcolor{gray} &   & \cellcolor{gray} & \cellcolor{gray} &  \\
    \hline
       $P_5$ & \cellcolor{gray} & \cellcolor{gray} &    & \cellcolor{gray} & \cellcolor{gray} &   & \cellcolor{gray} & \cellcolor{gray} & \cellcolor{gray} &   &   & \cellcolor{gray} & \cellcolor{gray} & \cellcolor{gray} & \cellcolor{gray}\\
    \hline
       $P_6$ & \cellcolor{gray} & \cellcolor{gray} & \cellcolor{gray} &   &   & \cellcolor{gray} & \cellcolor{gray} & \cellcolor{gray} &   & \cellcolor{gray} & \cellcolor{gray} & \cellcolor{gray} &   &   & \cellcolor{gray}\\
    \hline   
    \hline
       Date & 16 & 17 & 18 & 19 & 20 & 21 & 22 & 23 & 24 & 25 & 26 & 27 & 28 & 29 & 30\\
    \hline
       $P_4$ & \cellcolor{gray} & \cellcolor{gray} &   & \cellcolor{gray} & \cellcolor{gray} & \cellcolor{gray} &   &   & \cellcolor{gray} & \cellcolor{gray} & \cellcolor{gray} &   & \cellcolor{gray} & \cellcolor{gray} & \cellcolor{gray}\\
    \hline
       $P_5$ &   &   & \cellcolor{gray} & \cellcolor{gray} & \cellcolor{gray} &   & \cellcolor{gray} & \cellcolor{gray} &   & \cellcolor{gray} & \cellcolor{gray} & \cellcolor{gray} &   & \cellcolor{gray} & \cellcolor{gray}\\
    \hline
       $P_6$ & \cellcolor{gray} & \cellcolor{gray} & \cellcolor{gray} &   &   & \cellcolor{gray} & \cellcolor{gray} & \cellcolor{gray} & \cellcolor{gray} &   &   & \cellcolor{gray} & \cellcolor{gray} &   &  \\
    \hline   
    \end{tabular}}\\[1em]
    \item $U_3$: Day-shift group\\[1em]
    \resizebox{.9\textwidth}{!}{%
    \begin{tabular}{|c|c|c|c|c|c|c|c|c|c|c|c|c|c|c|c|}
    \hline
       Date & 1 & 2 & 3 & 4 & 5 & 6 & 7 & 8 & 9 & 10 & 11 & 12 & 13 & 14 & 15\\
    \hline
        $P_7$ &   &   &   &   & \cellcolor{gray} & \cellcolor{gray} & \cellcolor{gray} &   &   & \cellcolor{gray} &   & \cellcolor{gray} &   &   & \cellcolor{gray}\\
    \hline
        $P_8$ &   & \cellcolor{gray} & \cellcolor{gray} &   &   &   &   & \cellcolor{gray} & \cellcolor{gray} &   & \cellcolor{gray} & \cellcolor{gray} &   &   &  \\
    \hline
        $P_9$ &   & \cellcolor{gray} &   &   & \cellcolor{gray} & \cellcolor{gray} & \cellcolor{gray} &   & \cellcolor{gray} & \cellcolor{gray} & \cellcolor{gray} &   &   &   & \\
    \hline
        $P_{10}$ &   &   &   & \cellcolor{gray} & \cellcolor{gray} &   & \cellcolor{gray} &   &   &   &   &   & \cellcolor{gray} &   & \cellcolor{gray}\\
    \hline
        $P_{11}$ & \cellcolor{gray} &   &   & \cellcolor{gray} &   &   &   & \cellcolor{gray} &   &   & \cellcolor{gray} &   & \cellcolor{gray} & \cellcolor{gray} &  \\
    \hline
        $P_{12}$ & \cellcolor{gray} &   & \cellcolor{gray} &   &   & \cellcolor{gray} &   & \cellcolor{gray} & \cellcolor{gray} &   &   & \cellcolor{gray} &   & \cellcolor{gray} & \cellcolor{gray}\\
    \hline
        $P_{13}$ & \cellcolor{gray} & \cellcolor{gray} & \cellcolor{gray} & \cellcolor{gray} &   &   &   &   &   & \cellcolor{gray} &   &   & \cellcolor{gray} & \cellcolor{gray} &  \\
    \hline
    \hline
       Date & 16 & 17 & 18 & 19 & 20 & 21 & 22 & 23 & 24 & 25 & 26 & 27 & 28 & 29 & 30\\
    \hline
        $P_7$ &   & \cellcolor{gray} &   & \cellcolor{gray} &   & \cellcolor{gray} &   &   &   &   &   & \cellcolor{gray} & \cellcolor{gray} &   &  \\
    \hline
        $P_8$ & \cellcolor{gray} & \cellcolor{gray} & \cellcolor{gray} &   &   &   & \cellcolor{gray} & \cellcolor{gray} & \cellcolor{gray} & \cellcolor{gray} &   &   &   &   &  \\
    \hline
        $P_9$ &   &   &   & \cellcolor{gray} & \cellcolor{gray} & \cellcolor{gray} &   &   &   &   & \cellcolor{gray} &   & \cellcolor{gray} & \cellcolor{gray} & \cellcolor{gray}\\
    \hline
        $P_{10}$ &   & \cellcolor{gray} &   &   &   &   & \cellcolor{gray} &   & \cellcolor{gray} &   &   & \cellcolor{gray} &   &   &  \\
    \hline
        $P_{11}$ & \cellcolor{gray} &   & \cellcolor{gray} &   &   & \cellcolor{gray} & \cellcolor{gray} &   & \cellcolor{gray} & \cellcolor{gray} & \cellcolor{gray} & \cellcolor{gray} &   & \cellcolor{gray} & \cellcolor{gray}\\
    \hline
        $P_{12}$ & \cellcolor{gray} &   & \cellcolor{gray} & \cellcolor{gray} & \cellcolor{gray} &   &   & \cellcolor{gray} &   &   & \cellcolor{gray} &   &   &   &  \\
    \hline
        $P_{13}$ &   &   &   &   & \cellcolor{gray} &   &   & \cellcolor{gray} &   & \cellcolor{gray} &   &   & \cellcolor{gray} & \cellcolor{gray} & \cellcolor{gray}\\
    \hline
    \end{tabular}}
\end{itemize}

\item {Model 2: k = 5}
\begin{itemize}
    \item $U_1$: Graveyard-shift group\\[1em]
    \resizebox{.9\textwidth}{!}{%
    \begin{tabular}{|c|c|c|c|c|c|c|c|c|c|c|c|c|c|c|c|}
    \hline
       Date & 1 & 2 & 3 & 4 & 5 & 6 & 7 & 8 & 9 & 10 & 11 & 12 & 13 & 14 & 15\\
    \hline
        $P_1$ &   &   &   & \cellcolor{gray} & \cellcolor{gray} & \cellcolor{gray} & \cellcolor{gray} &   &   & \cellcolor{gray} & \cellcolor{gray} & \cellcolor{gray} & \cellcolor{gray} &   &   \\
    \hline
        $P_2$ & \cellcolor{gray} & \cellcolor{gray} & \cellcolor{gray} & \cellcolor{gray} & \cellcolor{gray} &   &   & \cellcolor{gray} & \cellcolor{gray} & \cellcolor{gray} & \cellcolor{gray} &   &   & \cellcolor{gray} & \cellcolor{gray}\\
    \hline
        $P_3$ & \cellcolor{gray} & \cellcolor{gray} & \cellcolor{gray} &   &   & \cellcolor{gray} & \cellcolor{gray} & \cellcolor{gray} & \cellcolor{gray} &   &   & \cellcolor{gray} & \cellcolor{gray} & \cellcolor{gray} & \cellcolor{gray}\\
    \hline
    \hline
        Date & 16 & 17 & 18 & 19 & 20 & 21 & 22 & 23 & 24 & 25 & 26 & 27 & 28 & 29 & 30\\
    \hline
        $P_1$ & \cellcolor{gray} & \cellcolor{gray} &   &   & \cellcolor{gray} & \cellcolor{gray} & \cellcolor{gray} & \cellcolor{gray} & \cellcolor{gray} &   & \cellcolor{gray} & \cellcolor{gray} & \cellcolor{gray} & \cellcolor{gray} & \cellcolor{gray}\\
    \hline
        $P_2$ & \cellcolor{gray} &   & \cellcolor{gray} & \cellcolor{gray} &   &   & \cellcolor{gray} & \cellcolor{gray} & \cellcolor{gray} & \cellcolor{gray} &   & \cellcolor{gray} & \cellcolor{gray} & \cellcolor{gray} & \cellcolor{gray}\\
    \hline
        $P_3$ &   & \cellcolor{gray} & \cellcolor{gray} & \cellcolor{gray} & \cellcolor{gray} & \cellcolor{gray} &   &   &   & \cellcolor{gray} & \cellcolor{gray} &   &   &   &   \\
    \hline
    \end{tabular}}\\[1em]
    \item $U_2$: Night-shift group\\[1em]
    \resizebox{.9\textwidth}{!}{%
    \begin{tabular}{|c|c|c|c|c|c|c|c|c|c|c|c|c|c|c|c|}
     \hline
       Date & 1 & 2 & 3 & 4 & 5 & 6 & 7 & 8 & 9 & 10 & 11 & 12 & 13 & 14 & 15\\
    \hline
        $P_4$ &   &   & \cellcolor{gray} & \cellcolor{gray} & \cellcolor{gray} & \cellcolor{gray} & \cellcolor{gray} &   &   & \cellcolor{gray} & \cellcolor{gray} &   & \cellcolor{gray} & \cellcolor{gray} &   \\
    \hline
        $P_5$ & \cellcolor{gray} & \cellcolor{gray} & \cellcolor{gray} & \cellcolor{gray} &   &   &   & \cellcolor{gray} & \cellcolor{gray} & \cellcolor{gray} & \cellcolor{gray} & \cellcolor{gray} &   &   & \cellcolor{gray}\\
    \hline
        $P_6$ & \cellcolor{gray} & \cellcolor{gray} &   &   & \cellcolor{gray} & \cellcolor{gray} & \cellcolor{gray} & \cellcolor{gray} & \cellcolor{gray} &   &   & \cellcolor{gray} & \cellcolor{gray} & \cellcolor{gray} & \cellcolor{gray}\\
    \hline
    \hline
        Date & 16 & 17 & 18 & 19 & 20 & 21 & 22 & 23 & 24 & 25 & 26 & 27 & 28 & 29 & 30\\
    \hline
        $P_4$ &   & \cellcolor{gray} & \cellcolor{gray} & \cellcolor{gray} & \cellcolor{gray} &   &   &   & \cellcolor{gray} & \cellcolor{gray} & \cellcolor{gray} &   & \cellcolor{gray} & \cellcolor{gray} & \cellcolor{gray}\\
    \hline
        $P_5$ & \cellcolor{gray} & \cellcolor{gray} & \cellcolor{gray} &   &   & \cellcolor{gray} & \cellcolor{gray} & \cellcolor{gray} & \cellcolor{gray} & \cellcolor{gray} &   & \cellcolor{gray} & \cellcolor{gray} &   &  \\
    \hline
        $P_6$ & \cellcolor{gray} &   &   & \cellcolor{gray} & \cellcolor{gray} & \cellcolor{gray} & \cellcolor{gray} & \cellcolor{gray} &   &   & \cellcolor{gray} & \cellcolor{gray} &   & \cellcolor{gray} & \cellcolor{gray} \\
    \hline
    \end{tabular}}\\[1em]
    \item $U_3$: Day-shift group\\[1em]
    \resizebox{.9\textwidth}{!}{%
    \begin{tabular}{|c|c|c|c|c|c|c|c|c|c|c|c|c|c|c|c|}
    \hline
       Date & 1 & 2 & 3 & 4 & 5 & 6 & 7 & 8 & 9 & 10 & 11 & 12 & 13 & 14 & 15\\
    \hline
        $P_7$ &   & \cellcolor{gray} &   &   & \cellcolor{gray} & \cellcolor{gray} &   &   &   & \cellcolor{gray} & 
  &   & \cellcolor{gray} &   &   \\
    \hline    
        $P_8$ &   &   & \cellcolor{gray} &   &   &   &   & \cellcolor{gray} &   &   &   & \cellcolor{gray} & \cellcolor{gray} &   &   \\
    \hline
        $P_9$ & \cellcolor{gray} & \cellcolor{gray} & \cellcolor{gray} & \cellcolor{gray} &   &   &   &   & \cellcolor{gray} &   & \cellcolor{gray} & \cellcolor{gray} &   &   & \cellcolor{gray}\\
    \hline
        $P_{10}$ & \cellcolor{gray} &   &   &   & \cellcolor{gray} &   & \cellcolor{gray} & \cellcolor{gray} & \cellcolor{gray} &   &   &   & \cellcolor{gray} & \cellcolor{gray} &  \\
    \hline
        $P_{11}$ & \cellcolor{gray} &   &   & \cellcolor{gray} &   & \cellcolor{gray} & \cellcolor{gray} & \cellcolor{gray} &   &   & \cellcolor{gray} & \cellcolor{gray} &   &   & \cellcolor{gray}\\
    \hline
        $P_{12}$ &   &   & \cellcolor{gray} &   &   &   & \cellcolor{gray} &   & \cellcolor{gray} & \cellcolor{gray} & \cellcolor{gray} &   &   & \cellcolor{gray} &  \\
    \hline
        $P_{13}$ &   & \cellcolor{gray} &   & \cellcolor{gray} & \cellcolor{gray} & \cellcolor{gray} &   &   &   & \cellcolor{gray} &   &   &   & \cellcolor{gray} & \cellcolor{gray} \\
    \hline
    \hline
        Date & 16 & 17 & 18 & 19 & 20 & 21 & 22 & 23 & 24 & 25 & 26 & 27 & 28 & 29 & 30\\
    \hline
        $P_7$ &   & \cellcolor{gray} &   & \cellcolor{gray} &   & \cellcolor{gray} &   &   & \cellcolor{gray} & \cellcolor{gray} &   & \cellcolor{gray} &   & \cellcolor{gray} & \cellcolor{gray}\\
    \hline
        $P_8$ &    &   &   & \cellcolor{gray} & \cellcolor{gray} &   &   &   &   &   &   &   &   &   &   \\
    \hline
        $P_9$ &   &   & \cellcolor{gray} & \cellcolor{gray} &   &   & \cellcolor{gray} &   &   & \cellcolor{gray} & \cellcolor{gray} &   & \cellcolor{gray} & \cellcolor{gray} &  \\
    \hline
        $P_{10}$ & \cellcolor{gray} &   & \cellcolor{gray} &   &   &   & \cellcolor{gray} & \cellcolor{gray} &   &  &   & \cellcolor{gray} & \cellcolor{gray} & \cellcolor{gray} &  \\
    \hline
        $P_{11}$ & \cellcolor{gray} & \cellcolor{gray} &   &   & \cellcolor{gray} & \cellcolor{gray} & \cellcolor{gray} &   &   &   & \cellcolor{gray} &   &   &   & \cellcolor{gray} \\
    \hline
        $P_{12}$ & \cellcolor{gray} &   &   &   &   & \cellcolor{gray} &   & \cellcolor{gray} & \cellcolor{gray} &   & \cellcolor{gray} & \cellcolor{gray} &   &   & \cellcolor{gray} \\
    \hline
        $P_{13}$ &   & \cellcolor{gray} & \cellcolor{gray} &   & \cellcolor{gray} &   &   & \cellcolor{gray} & \cellcolor{gray} & \cellcolor{gray} &   &   & \cellcolor{gray} &   &   \\
    \hline
    \end{tabular}}\\[1em]
\end{itemize}
\end{itemize}
\end{example}

Based on our analysis for both models 1 and 2, we have observed that all the rules mentioned in Section 3 are satisfied. This confirms that using the QUBO framework is an effective way to solve the NSP. Moreover, despite having different designs in terms of consecutive working days, model 1 and model 2 generate considerably different shift tables. This indicates that the NSP is highly sensitive to changes. If labor laws are changed in the real world without an efficient algorithm for scheduling, companies and management divisions may face difficulties. Hence, we demonstrate that our method can be applied to handle various scenarios, 
where our mathematical model created in Section 3  proposes a general form for solving the NSP. This general framework keeps its flexibility and feasibility to fit the needs of different hospitals by suitably adjusting parameters and adding additional constraints. 
Additionally, this achievement can establish a benchmark for using the quantum annealing algorithm to solve the NSP by applying the same QUBO formulation.

\begin{example}
%\subsection{Discussion}
In our previous example, we only consider hard constraints that must be satisfied. To improve job satisfaction, we want to maximize the number of two-day leaves for the nurses in graveyard shifts by adding one more soft constraint to our objective function 
\begin{equation} 
\label{concentrateleave}
    H_{cost}(\mathbf{x}) = \sum_{i=1}^{m_1} \sum_{j = 1}^{d-1} 1 - x_{ij}x_{i,j+1}.
\end{equation}
This newly introduced constraint is a common preference among all nurses, as consecutive days off allow them to plan longer-term activities for relaxation. However, this requirement is not essential from the hospital's perspective, so it is classified as a soft constraint. The design for two-day leaves utilizes a multiplication approach, where the product will be zero if one of the binary variables is zero. Additionally, the objective is to maximize the number of two-day leaves, effectively minimizing the sum of $1 - x_{i, j} x_{i, j+1}$. This formula encourages more consecutive days off in the scheduling algorithm.
The following tables show the output results concerning different $k$ values.

% Considering that for $i = 1$, (\ref{concentrateleave}) becomes 
% \begin{equation}
%     H_{cost}(\{x_{1,j}\}) = \sum_{j = 1}^{d-1} 1 - x_{ij}x_{1,j+1} = (d-1) - \sum_{j = 1}^{d-1} x_{1,j}x_{1,j+1}
% \end{equation}
%Noticed that the total number of $0$ is fixed. This formula gets lower energy when $0$'s are gathered. 
\begin{itemize}
    \item Model 1: $k = 4$\\ Graveyard-shift group\\
    \resizebox{.9\textwidth}{!}{%
    \begin{tabular}{|c|c|c|c|c|c|c|c|c|c|c|c|c|c|c|c|}
    \hline
       Date & 1 & 2 & 3 & 4 & 5 & 6 & 7 & 8 & 9 & 10 & 11 & 12 & 13 & 14 & 15\\
    \hline
       $P_1$ &   &   & \cellcolor{gray} & \cellcolor{gray} &   &   & \cellcolor{gray} & \cellcolor{gray} & \cellcolor{gray} &   & \cellcolor{gray} & \cellcolor{gray} & \cellcolor{gray} & \cellcolor{gray} &  \\
    \hline
       $P_2$ & \cellcolor{gray} & \cellcolor{gray} &   & \cellcolor{gray} & \cellcolor{gray} & \cellcolor{gray} &   & \cellcolor{gray} & \cellcolor{gray} & \cellcolor{gray} &   &   & \cellcolor{gray} & \cellcolor{gray} & \cellcolor{gray} \\
    \hline
       $P_3$ & \cellcolor{gray} & \cellcolor{gray} & \cellcolor{gray} &   & \cellcolor{gray} & \cellcolor{gray} & \cellcolor{gray} &  &   & \cellcolor{gray} & \cellcolor{gray} & \cellcolor{gray} &  &   & \cellcolor{gray}\\
    \hline   
    \hline
       Date & 16 & 17 & 18 & 19 & 20 & 21 & 22 & 23 & 24 & 25 & 26 & 27 & 28 & 29 & 30\\
    \hline
       $P_1$ &   & \cellcolor{gray} & \cellcolor{gray} & \cellcolor{gray} & \cellcolor{gray} &   &   & \cellcolor{gray} & \cellcolor{gray} & \cellcolor{gray} & \cellcolor{gray} &   &   & \cellcolor{gray} & \cellcolor{gray}\\
    \hline
       $P_2$ & \cellcolor{gray} &   &   & \cellcolor{gray} & \cellcolor{gray} & \cellcolor{gray} & \cellcolor{gray} &   &   &  \cellcolor{gray} & \cellcolor{gray} & \cellcolor{gray} & \cellcolor{gray} &   &  \\
    \hline
       $P_3$ & \cellcolor{gray} & \cellcolor{gray} & \cellcolor{gray} &   &   & \cellcolor{gray} & \cellcolor{gray} & \cellcolor{gray} & \cellcolor{gray} &   &   & \cellcolor{gray} & \cellcolor{gray} & \cellcolor{gray} & \cellcolor{gray}\\
    \hline   
    \end{tabular}}\\[1em]
    \item Model 2: $k = 5$ \\ Graveyard-shift group\\
    \resizebox{.9\textwidth}{!}{%
    \begin{tabular}{|c|c|c|c|c|c|c|c|c|c|c|c|c|c|c|c|}
    \hline
       Date & 1 & 2 & 3 & 4 & 5 & 6 & 7 & 8 & 9 & 10 & 11 & 12 & 13 & 14 & 15\\
    \hline
       $P_1$ & \cellcolor{gray} & \cellcolor{gray} & \cellcolor{gray} &   &   & \cellcolor{gray} & \cellcolor{gray} & \cellcolor{gray} & \cellcolor{gray} & \cellcolor{gray} &   &  &  \cellcolor{gray} & \cellcolor{gray} & \cellcolor{gray}\\
    \hline
       $P_2$ &   &   &   & \cellcolor{gray} & \cellcolor{gray} & \cellcolor{gray} & \cellcolor{gray} & \cellcolor{gray} &  &   & \cellcolor{gray} & \cellcolor{gray} & \cellcolor{gray} & \cellcolor{gray} & \cellcolor{gray}\\
    \hline
       $P_3$ & \cellcolor{gray} & \cellcolor{gray} & \cellcolor{gray} & \cellcolor{gray} & \cellcolor{gray} &   &   &   & \cellcolor{gray} & \cellcolor{gray} & \cellcolor{gray} & \cellcolor{gray} &   &   &  \\
    \hline   
    \hline
       Date & 16 & 17 & 18 & 19 & 20 & 21 & 22 & 23 & 24 & 25 & 26 & 27 & 28 & 29 & 30\\
    \hline
       $P_1$ & \cellcolor{gray} &   &   &   & \cellcolor{gray} & \cellcolor{gray} & \cellcolor{gray} & \cellcolor{gray} & \cellcolor{gray} &   &   & \cellcolor{gray} & \cellcolor{gray} & \cellcolor{gray} & \cellcolor{gray}\\
    \hline
       $P_2$ &   & \cellcolor{gray} & \cellcolor{gray} & \cellcolor{gray} & \cellcolor{gray} & \cellcolor{gray} &   &   &   & \cellcolor{gray} & \cellcolor{gray} & \cellcolor{gray} & \cellcolor{gray} &   &  \\
    \hline
       $P_3$ & \cellcolor{gray} & \cellcolor{gray} & \cellcolor{gray} & \cellcolor{gray} &   &   & \cellcolor{gray} & \cellcolor{gray} & \cellcolor{gray} & \cellcolor{gray} & \cellcolor{gray} &   &   & \cellcolor{gray} & \cellcolor{gray}\\
    \hline   
    \end{tabular}}\\[1em]
\end{itemize}
\end{example}

 Since the number of nurses is predetermined and the constraints are satisfied in our calculation, we can see that the total number of white blocks representing the leaves of nurses should be fixed. To minimize~\eqref{concentrateleave}, we should increase the number of consecutive working days. In this case, the number of consecutive leaves for the nurses will increase simultaneously, as seen in our above two results.
From the outcomes presented in Models 1 and 2 for Examples 1 and 2, we observe an increase in two-day leaves, aligning with our expectations. This observation demonstrates the efficacy and adaptability of incorporating soft constraints into our model to meet the needs of nurses. While the working days for each nurse remain similar, the addition of consecutive day leaves in this example highlights the model's capability to adapt with minimal modifications. This approach effectively addresses specific scenarios once the schedule framework is established.

Real-world scheduling involves many diverse requirements, and the given example highlights the capacity of the QUBO framework to incorporate an additional requirement without violating fundamental rules. When it comes to practical settings, nurses often have different shift preferences. Employing this efficient technique for generating shift schedules can significantly reduce the workload of scheduling staff, streamline the process, and ensure that nurse preferences are adequately considered.

\section{Conclusion}\label{sec6}
We, in this work, propose a simple and readily implemented computing model to solve the nurse scheduling problem. We first reduce the number of variables by separating nurses into three groups. We then solve this model by using a built-in solver, the SAA, in PyQUBO. We show that this model can not only be solved easily but can also be modified to satisfy additional constraints without difficulty. Though the initial group separation can reduce the number of variables, a clear allocation of nursing resources into three groups might not exist in some cases. For example, a nurse might prefer graveyard shifts but can also support day or night shifts. To this, we can modify our proposed model to present diverse conditions by adding more variables. However, while facing various scenarios, solving this problem efficiently without increasing the number of variables is worthy of further study.

%%  The bibliography

%\begin{thebibliography}{9}
%%  Use \bibitem{r1} or \bibitem[Surname(2010)]{r1} (for authoryear case)
\bibliographystyle{imsart-nameyear}
\bibliography{bibliography}

%\end{thebibliography}

\end{document}